\documentstyle[amssymb,12pt]{article}

\newtheorem{th}{Theorem}[section]

\newtheorem{prop}[th]{Proposition}
\newtheorem{cor}[th]{Corollary}
\newtheorem{defn}[th]{Definition}
\newenvironment{defn-new}{\begin{defn} \em}{\end{defn}}
\newtheorem{rem}[th]{Remark}
\newenvironment{rem-new}{\begin{rem} \em}{\end{rem}}
\newtheorem{ex}[th]{Example}
\newenvironment{ex-new}{\begin{ex} \em}{\end{ex}}
\newtheorem{prob}[th]{Problem}
\newenvironment{prob-new}{\begin{prob} \em}{\end{prob}}

\newenvironment{notation-new}{\begin{rem} \em}{\end{rem}}

\newenvironment{agr-new}{\begin{rem} \em}{\end{rem}}

\makeatletter \@addtoreset{equation}{section} \makeatother
\begin{document}

\begin{center}
{\Large {\bf Geometry of lightlike hypersurfaces of a statistical manifold}}%
\bigskip \bigskip

O\u{g}uzhan Bahad\i r, Mukut Mani Tripathi
\end{center}

\bigskip \bigskip



\medskip

\noindent {\bf Abstract.} Lightlike hypersurfaces of a statistical manifold
are studied. It is shown that a lightlike hypersurface of a statistical
manifold is not a statistical manifold with respect to the induced
connections, but the screen distribution has a canonical statistical
structure. Some relations between induced geometric objects with respect to
dual connections in a lightlike hypersurface of a statistical manifold are
obtained. An example is presented. Induced Ricci tensors for lightlike
hypersurface of a statistical manifold are computed.

\section{Introduction\label{sect-intro}}

A statistical manifold, the Riemannian connection used to model the
information, the fields of information geometry, as such a generalization of
the Riemannian manifold equipped with a relatively new mathematics branch,
uses the differential geometry tool to examine the statistical inference,
information loss and prediction \cite{Calin-Udr-2014-book}. In 1975, The
role of differential geometry in statistics was first emphasized by Efron
\cite{Efron-1975}. Later, Amari used differential geometric tools to develop
this idea \cite{Amari-1982}, \cite{Amari-1985-book}.

In 1989, Vos \cite{Vos-1989} initiated the study of geometry of submanifolds
of statistical manifolds. He obtained Gauss-Weingarten formulas, Gauss and
Codazzi equations, etc.. Later, in 2009, Furuhata \cite{Furuhata-2009}
studied hypersurfaces of a statistical manifold. Also, Aydin et. al. studied
submanifolds of statistical manifolds of constant curvature \cite%
{Aydin-MM-2015}.

On the other hand, lightlike geometry is one of the important research areas
in differential geometry and has many applications in physics and
mathematics. The geometry of lightlike submanifolds of a semi-Riemannian
manifold was presented by K.L. Duggal and A. Bejancu in \cite%
{Duggal-Bejancu-1996} (see also \cite{Duggal-Jin-2007}, \cite%
{Duggal-Sahin-2010}). Lightlike hypersurfaces in various spaces have been
studied by many authors including those of \cite{Atindogbe-ET-2006}, \cite%
{Duggal-2012}, \cite{Duggal-Jin-2007}, \cite{Izumiya-2015}, \cite%
{Kazan-Sah-2016}, \cite{Kilic-Bah-2012}, \cite{Liu-Pan-2015}, \cite%
{Massamba-2008-DGDS}, \cite{Massamba-2008-Turk}, \cite{Massamba-2013}, \cite%
{Poy-Yasar-2017}.

Motivated by these circumstances, in this paper, we initiate the study of
lightlike geometry of statistical manifolds. In section~\ref{sect-prel}, we
present basic definitions and results about statistical manifolds and
lightlike hypersurfaces. In Section~\ref{sect-lightlike}, we show that
induced connections on a lightlike hypersurface of a statistical manifold
are not dual connections and a lightlike hypersurface is not statistical
manifold. Moreover, we show that the second fundamental forms are not
degenerate. Later, we characterize the parallelness and integrability of the
screen distribution. Equivalent conditions are also obtained between the
induced objects. This section concludes with an example. In section~\ref%
{sect-curvature}, we obtain formula for curvature tensors of a lightlike
hypersurface of a statistical manifold. In general, in lightlike geometry,
Ricci tensor is not symmetric, so we also obtain new conditions for Ricci
tensor to be symmetric.

\section{Preliminaries\label{sect-prel}}

We begin with the following definition.

\begin{defn-new}
{\rm \cite{Furuhata-2009}} Let $\widetilde{M}$ be a smooth manifold. Let $%
\widetilde{D}$ be an affine connection with the torsion tensor $T^{%
\widetilde{D}}$ and $\widetilde{g}$ a semi-Riemannian metric on $\widetilde{M%
}$. Then the pair $(\widetilde{D},\widetilde{g})$ is called a statistical
structure on $\widetilde{M}$ if
\begin{enumerate}
\item[{\rm (1)}] $(\widetilde{D}_{X}\widetilde{g})(Y,Z) - (\widetilde{D}_{Y}%
\widetilde{g})(X,Z) = \widetilde{g}(T^{\widetilde{D}}(X,Y),Z)$ \newline
for all $X,Y,Z\in \Gamma(T\widetilde{M})$, and

\item[{\rm (2)}] $T^{\widetilde{D}}=0$.
\end{enumerate}
\end{defn-new}

\begin{defn-new}
Let $(\widetilde{M},\widetilde{g})$ be a semi-Riemannian manifold. Two
affine connections $\widetilde{D}$ and $\widetilde{D}^{\ast }$ on $%
\widetilde{M}$ are said to be dual with respect to the metric $\widetilde{g}$%
, if
\begin{equation}
Z\widetilde{g}(X,Y)=\widetilde{g}(\widetilde{D}_{Z}X,Y)+\widetilde{g}(X,%
\widetilde{D}_{Z}^{\ast }Y)  \label{eq-dual-con}
\end{equation}%
for all $X,Y,Z\in \Gamma (T\widetilde{M})$.
\end{defn-new}

A statistical manifold will be represented by $(\widetilde{M},\widetilde{g},%
\widetilde{D},\widetilde{D}^{\ast })$. If $\widetilde{D}^{0}$ is Levi-Civita
connection of $\widetilde{g}$, then
\begin{equation}
\widetilde{D}^{0}=\frac{1}{2}(\widetilde{D}+\widetilde{D}^{\ast }).
\label{le}
\end{equation}%
In (\ref{eq-dual-con}), if we choose $\widetilde{D}^{\ast }=\widetilde{D}$
then Levi-Civita connection is obtained.

\medskip

Let $(M,g)$ be a submanifold of $(\widetilde{M},\widetilde{g})$. If $%
(M,g,D,D^{\ast })$ is a statistical manifold, then $(M,g,D,D^{\ast })$ is
called a statistical submanifold of $(\widetilde{M},\widetilde{g},\widetilde{%
D},\widetilde{D}^{\ast })$, where $D$, $D^{\ast }$ are affine dual
connections on $M$ and $\widetilde{D}$, $\widetilde{D}^{\ast }$ are affine
dual connections on $\widetilde{M}$ (see \cite{Amari-1985-book}, \cite%
{Furuhata-2009},\cite{Vos-1989}).

\medskip

Now, let $(\bar{M},\bar{g})$ be an $(m+2)$-dimensional semi-Riemannian
manifold with ${\rm index}(\bar{g})=q\geq 1$. Let $(M,g)$ be a hypersurface
of $(\bar{M},\bar{g})$ with $g=\bar{g}|_{M}$. If the induced metric $g$ on $%
M $ is degenerate, then $M$ is called a lightlike (null or degenerate)
hypersurface (\cite{Duggal-Bejancu-1996}, \cite{Duggal-Jin-2007}, \cite%
{Duggal-Sahin-2010}). In this case, there exists a null vector field $\xi
\neq 0$ on $M$ such that
\begin{equation}
g\left( \xi ,X\right) =0,\qquad \forall \;X\in \Gamma \left( TM\right) .
\label{eq-null-1}
\end{equation}%
The radical or the null space of $T_{x}M$, at each point $x\in M$, is a
subspace $Rad~T_{x}M$ defined by
\begin{equation}
Rad~T_{x}M=\{\xi \in T_{x}M:g_{x}(\xi ,X)=0,\;X\in \Gamma (TM)\}.
\label{eq-null-2}
\end{equation}%
The dimension of $Rad~T_{x}M$ is called the nullity degree of $g$. We recall
that the nullity degree of $g$ for a lightlike hypersurface of $(\bar{M},%
\bar{g})$ is $1$. Since $g$ is degenerate and any null vector being
orthogonal to itself, $T_{x}M^{\perp }$ is also null and
\begin{equation}
Rad~T_{x}M=T_{x}M\cap T_{x}M^{\perp }.  \label{eq-null-3}
\end{equation}%
Since $\dim T_{x}M^{\perp }=1$ and $\dim Rad~T_{x}M=1,$ we have $%
Rad~T_{x}M=T_{x}M^{\perp }$. We call $Rad~TM$ a radical distribution and it
is spanned by the null vector field $\xi $. The complementary vector bundle $%
S(TM)$ of $Rad~TM$ in $TM$ is called the screen bundle of $M$. We note that
any screen bundle is non-degenerate. This means that
\begin{equation}
TM=Rad~TM\perp S(TM),  \label{eq-null-4}
\end{equation}%
with $\perp $ denoting the orthogonal-direct sum. The complementary vector
bundle $S(TM)^{\perp }$ of $S(TM)$ in $T\bar{M}$ is called screen
transversal bundle and it has rank $2$. Since $Rad~TM$ is a lightlike
subbundle of $S(TM)^{\perp }$ there exists a unique local section $N$ of $%
S(TM)^{\perp }$ such that
\begin{equation}
\bar{g}(N,N)=0,\quad \bar{g}(\xi ,N)=1.  \label{eq-null-5}
\end{equation}%
Note that $N$ is transversal to $M$ and $\{\xi ,N\}$ is a local frame field
of $S(TM)^{\perp }$ and there exists a line subbundle $ltr(TM)$ of $T\bar{M}$
, and it is called the lightlike transversal bundle, locally spanned by $N$.
Hence we have the following decomposition:
\begin{equation}
T\bar{M}=TM\oplus ltr(TM)=S(TM)\bot Rad~TM\oplus ltr(TM),  \label{eq-null-6}
\end{equation}%
where $\oplus $ is the direct sum but not orthogonal (\cite%
{Duggal-Bejancu-1996}, \cite{Duggal-Jin-2007}). From the above decomposition
of a semi-Riemannian manifold $\bar{M}$ along a lightlike hypersurface $M$,
we can consider the local quasi-orthonormal field of frames of $\bar{M}$
along $M$ given by
\[
\{E_{1},\ldots ,E_{m},\xi ,N\},
\]%
where $\{E_{1},\ldots ,E_{m}\}$ is an orthonormal basis of $\Gamma (S(TM))$.
In view of the splitting $(\ref{eq-null-6})$, we have the following Gauss
and Weingarten formulas, respectively,
\begin{equation}
\bar{\nabla}_{X}Y=\nabla _{X}Y+h(X,Y),  \label{eq-null-7}
\end{equation}%
\begin{equation}
\bar{\nabla}_{X}N=-A_{N}X+\nabla _{X}^{t}N  \label{eq-null-8}
\end{equation}%
for any $X,Y\in \Gamma (TM)$, where $\nabla _{X}Y,\ A_{N}X\in \Gamma (TM)$
and $h(X,Y),\ \nabla _{X}^{t}N\in \Gamma (ltr(TM))$. If we set
\[
B(X,Y)=\bar{g}(h(X,Y),\xi )\quad {\rm and}\quad \tau (X)=\bar{g}(\nabla
_{X}^{t}N,\xi ),
\]%
then (\ref{eq-null-7}) and (\ref{eq-null-8}) become
\begin{equation}
\overline{\nabla }_{X}Y=\nabla _{X}Y+B(X,Y)N,  \label{eq-null-9}
\end{equation}%
\begin{equation}
\overline{\nabla }_{X}N=-A_{N}X+\tau (X)N,  \label{eq-null-10}
\end{equation}%
respectively. Here, $B$ and $A$ are called the second fundamental form and
the shape operator of the lightlike hypersurface $M$, respectively \cite%
{Duggal-Bejancu-1996}. Let $P$ be the projection of $S(TM)$ on $M$. Then,
for any $X\in \Gamma (TM)$, we can write
\begin{equation}
X=PX+\eta (X)\xi ,  \label{eq-null-11}
\end{equation}%
where $\eta $ is a $1$-form given by
\begin{equation}
\eta (X)=\bar{g}(X,N).  \label{eq-null-12}
\end{equation}%
From (\ref{eq-null-9}), we have
\begin{equation}
(\nabla _{X}g)(Y,Z)=B(X,Y)\eta (Z)+B(X,Z)\eta (Y),  \label{eq-null-13}
\end{equation}%
for all $X,Y,Z\in \Gamma (TM)$, where the induced connection $\nabla $ is a
non-metric connection on $M$. From (\ref{eq-null-4}), we have
\begin{equation}
\nabla _{X}W=\nabla _{X}^{\ast }W+h^{\ast }(X,W)=\nabla _{X}^{\ast
}W+C(X,W)\xi ,  \label{eq-null-14}
\end{equation}%
\begin{equation}
\nabla _{X}\xi =-A_{\xi }^{\ast }X-\tau (X)\xi  \label{eq-null-15}
\end{equation}%
for all $X\in \Gamma (TM)$, $W\in \Gamma (S(TM))$, where $\nabla _{X}^{\ast
}W$ and $A_{\xi }^{\ast }X$ belong to $\Gamma (S(TM))$. Here $C$, $A_{\xi
}^{\ast }$ and $\nabla ^{\ast }$ are called the local second fundamental
form, the local shape operator and the induced connection on $S(TM)$,
respectively. We also have
\begin{equation}
g(A_{\xi }^{\ast }X,W)=B(X,W),\ g(A_{\xi }^{\ast }X,N)=0,\ B(X,\xi )=0,\ \ \
g(A_{N}X,N)=0.  \label{eq-null-16}
\end{equation}%
Moreover, from the first and third equations of (\ref{eq-null-16}), we have
\begin{equation}
A_{\xi }^{\ast }\xi =0.  \label{eq-null-17}
\end{equation}%
The mean curvature $H$ of $M$ with respect to an $\{E_{i}\},\ i=1,\ldots m,$
orthonormal basis of $\Gamma (S(TM))$ is defined by
\begin{equation}
H=\frac{1}{m}\sum_{i=1}^{m}\varepsilon _{i}B(E_{i},E_{i}),\quad \varepsilon
_{i}=g(E_{i},E_{i}).  \label{eq-null-18-H}
\end{equation}%
Let $x\in M$ and $\Pi = {\rm span}\{E_{i},E_{j}\}$ be a $2$-dimensional
non-degenerate plane of $T_{x}M$. The sectional curvature of $\Pi$ at $x\in
M $ is defined by \cite{Beem-EE-1996-book}
\begin{equation}
\kappa _{ij}=\frac{g(R(E_{j},E_{i})E_{j},E_{i})}{%
g(E_{i},E_{i})g(E_{j},E_{j})-g(E_{j},E_{i})^{2}},  \label{eq-kes-1}
\end{equation}
Now, let $x\in M$ and $\xi$ be a null vector of $T_{x}M$. A plane $\Pi$ of $%
T_{x}M$ is a null plane if it contains $\xi$ and $E_{i}$ such that $\bar{g}%
(\xi,E_{i})=0$ and $g(E_{i},E_{i})=\varepsilon_{i}$. Then the null sectional
curvature is given by \cite{Beem-EE-1996-book}
\begin{equation}
\kappa_{i}^{null} = \frac{g(R_{u}(\xi,E_{i})\xi,E_{i})} {g_{u}(E_{i},E_{i})}.
\label{eq-kes-2}
\end{equation}

\section{Lightlike hypersurfaces of a statistical manifold\label%
{sect-lightlike}}

Let $(M,g)$ be a lightlike hypersurface of a statistical manifold $(%
\widetilde{M},\widetilde{g},\widetilde{D},\widetilde{D}^{\ast })$. Then,
Gauss and Weingarten formulas with respect to dual connections are given by
\cite{Furuhata-2009}
\begin{equation}
\widetilde{D}_{X}Y=D_{X}Y+B(X,Y)N,  \label{eq-D-tilde-X-Y}
\end{equation}%
\begin{equation}
\widetilde{D}_{X}N=-A_{N}^{\ast }X+\tau ^{\ast }(X)N,  \label{eq-D-tilde-X-N}
\end{equation}%
\begin{equation}
\widetilde{D}_{X}^{\ast }Y=D_{X}^{\ast }Y+B^{\ast }(X,Y)N,
\label{eq-D-tilde-star-X-Y}
\end{equation}%
\begin{equation}
\widetilde{D}_{X}^{\ast }N=-A_{N}X+\tau (X)N  \label{eq-D-tilde-star-X-N}
\end{equation}%
for all $X,Y\in \Gamma (TM),\;N\in \Gamma (ltrTM)$, where $D_{X}Y$,$\
D_{X}^{\ast }Y$, $A_{N}X$,$\ A_{N}^{\ast }X\in \Gamma (TM)$ and
\[
B(X,Y)=\widetilde{g}(\widetilde{D}_{X}Y,\xi ),\quad {\tau }^{\ast }(X)=%
\widetilde{g}(\widetilde{D}_{X}N,\xi ),
\]%
\[
B^{\ast }(X,Y)=\widetilde{g}(\widetilde{D}_{X}^{\ast }Y,\xi ),\quad {\tau }%
(X)=\widetilde{g}(\widetilde{D}_{X}^{\ast }N,\xi ).
\]%
Here, $D$, $D^{\ast }$, $B$, $B^{\ast }$, ${A}_{N}$ and $A_{N}^{\ast }$ are
called the induced connections on $M$, the second fundamental forms and the
Weingarten mappings with respect to $\widetilde{D}$ and $\widetilde{D}^{\ast
}$, respectively. Using Gauss formulas and the equation (\ref{eq-dual-con}),
we obtain
\begin{eqnarray}
Xg(Y,Z) &=&g(\widetilde{D}_{X}Y,Z)+g(Y,\widetilde{D}_{X}^{\ast }Z),\
\nonumber \\
&=&g(D_{X}Y,Z)+g(Y,D_{X}^{\ast }Z)+B(X,Y)\eta (Z)+B^{\ast }(X,Z)\eta (Y).
\label{eq-LH-SM-1}
\end{eqnarray}

From the equation (\ref{eq-LH-SM-1}), we have the following result.

\begin{th}
\label{th-LH-SM-1} Let $(M,g)$ be a lightlike hypersurface of a statistical
manifold $(\widetilde{M},\widetilde{g},\widetilde{D},\widetilde{D}^{\ast })$%
. Then\/{\rm :}

\begin{enumerate}
\item[{\bf (i)}] Induced connections $D$ and $D^{\ast }$ are not dual
connections.

\item[{\bf (ii)}] A lightlike hypersurface of a statistical manifold need
not a statistical manifold with respect to the dual connections.
\end{enumerate}
\end{th}

Using Gauss and Weingarten formulas in (\ref{eq-LH-SM-1}), we get
\begin{eqnarray}
(D_{X}g)(Y,Z) + (D_{X}^{\ast }g)(Y,Z) &=& B(X,Y)\eta (Z) + B(X,Z)\eta(Y)
\nonumber \\
&&+B^{\ast }(X,Y)\eta (Z)+B^{\ast }(X,Z)\eta (Y).  \label{th-LH-SM-2}
\end{eqnarray}

\begin{prop}
\label{pro5} Let $(M,g)$ be a lightlike hypersurface of a statistical
manifold $(\widetilde{M},\widetilde{g},\widetilde{D},\widetilde{D}^{\ast })$%
. Then the following assertions are true\/{\rm :}

\begin{enumerate}
\item[{\bf (i)}] Induced connections $D$ and $D^{\ast }$ are symmetric
connection.

\item[{\bf (ii)}] The second fundamental forms $B$ and $B^{\ast }$ are
symmetric.
\end{enumerate}
\end{prop}

\noindent {\bf Proof.} We know that $T^{\widetilde{D}}=0$. Moreover,
\begin{eqnarray}
T^{\widetilde{D}}(X,Y)&=&\widetilde{D}_{X}Y-\widetilde{D}_{Y}X-[X,Y]
\nonumber \\
&=&D_{X}Y-D_{Y}X-[X,Y]+B(X,Y)N-B(Y,X)N=0.  \label{th-LH-SM-3}
\end{eqnarray}
Comparing the tangent and transversal components of (\ref{th-LH-SM-3}), we
obtain
\[
B(X,Y)=B(Y,X), \qquad T^{D}=0,
\]
where $T^{D}$ is the torsion tensor field of $D$. Thus, second fundamental
form $B$ is symmetric and induced connection $D$ is symmetric connection.

Similarly, it can be shown that the second fundamental form $B^{\ast}$ is
symmetric and the induced connection $D^{\ast}$ is a symmetric connection. $%
\blacksquare$

\medskip

Let $P$ denote the projection morphism of $\Gamma (TM)$ on $\Gamma (S(TM))$
with respect to the decomposition (\ref{eq-null-4}). Then, we have
\begin{equation}
D_{X}PY=\nabla _{X}PY+\overline{h}(X,PY),
\end{equation}%
\begin{equation}
D_{X}\xi =-\overline{A}_{\xi }X+\overline{\nabla }_{X}^{t}\xi =0
\end{equation}%
for all $X,Y\in \Gamma (TM)$ and $\xi \in \Gamma (RadTM)$, where $\nabla
_{X}PY$ and $\overline{A}_{\xi }X$ belong to $\Gamma (S(TM))$, $\nabla $ and
$\overline{\nabla }^{t}$ are linear connections on $\Gamma (S(TM))$ and $%
\Gamma (RadTM)$ respectively. Here, $\overline{h}$ and $\overline{A}$ are
called screen second fundamental form and screen shape operator of $S(TM)$,
respectively. If we define
\begin{equation}
C(X,PY)=g(\overline{h}(X,PY),N),
\end{equation}%
\begin{equation}
\varepsilon (X)=g(\overline{\nabla }_{X}^{t}\xi ,N),\;\forall X,Y\in \Gamma
(TM).
\end{equation}%
One can show that
\[
\varepsilon (X)=-\tau (X).
\]%
Therefore, we have
\begin{equation}
D_{X}PY=\nabla _{X}PY+C(X,PY)\xi ,  \label{c}
\end{equation}%
\begin{equation}
D_{X}\xi =-\overline{A}_{\xi }X-\tau (X)\xi =0,\;\forall X,Y\in \Gamma (TM).
\label{31}
\end{equation}%
Here $C(X,PY)$ is called the local screen fundamental form of $S(TM)$.

\medskip

Similarly, the relations of induced dual objects on $S(TM)$ are given by
\begin{equation}
D_{X}^{\ast }PY=\nabla _{X}^{\ast }PY+C^{\ast }(X,PY)\xi ,  \label{d}
\end{equation}%
\begin{equation}
D_{X}^{\ast }\xi =-\overline{A}_{\xi }^{\ast }X-\tau ^{\ast }(X)\xi
=0,\;\forall X,Y\in \Gamma (TM).  \label{33}
\end{equation}%
Using (\ref{eq-LH-SM-1}), (\ref{c}), (\ref{d}) and Gauss-Weingarten
formulas, the relationship between induced geometric objects are given by
\begin{equation}
B(X,\xi )+B^{\ast }(X,\xi )=0,\;g(A_{N}X+A_{N}^{\ast }X,N)=0,  \label{34}
\end{equation}%
\begin{equation}
C(X,PY)=g(A_{N}X,PY),\;C^{\ast }(X,PY)=g(A_{N}^{\ast }X,PY).  \label{18}
\end{equation}

Now, using the equation (\ref{34}) we can state the following result.

\begin{prop}
\label{pr6} Let $(M,g)$ be a lightlike hypersurface of a statistical
manifold $(\widetilde{M},\widetilde{g},\widetilde{D},\widetilde{D}^{\ast })$%
. Then second fundamental forms $B$ and $B^{\ast }$ are not degenerate.
\end{prop}

Additionally, due to $\widetilde{D}$ and $\widetilde{D}^{\ast }$ are dual
connections we obtain
\begin{equation}
B(X,Y)=g(\overline{A}_{\xi }^{\ast }X,Y)+B^{\ast }(X,\xi ),  \label{7}
\end{equation}%
\begin{equation}
B^{\ast }(X,Y)=g(\overline{A}_{\xi }X,Y)+B(X,\xi ).  \label{8}
\end{equation}%
Using (\ref{7}) and (\ref{8}) we get
\[
\overline{A}_{\xi }^{\ast }\xi +\overline{A}_{\xi }\xi =0.
\]

\begin{prop}
\label{pr7} Let $(M,g)$ be a lightlike hypersurface of a statistical
manifold $(\widetilde{M},\widetilde{g},\widetilde{D},\widetilde{D}^{\ast })$%
. Then the screen distribution $(S(TM),g,\nabla,\nabla^{\ast })$ has a
statistical structure.
\end{prop}

\noindent {\bf Proof.} From (\ref{eq-LH-SM-1}), for any $X,Y\in \Gamma
(S(TM))$ we obtain
\[
Xg(Y,Z)=g(D_{X}Y,Z)+g(Y,D_{X}^{\ast }Z).
\]%
Using (\ref{d}) in the last equation, we get
\[
Xg(Y,Z)=g(\nabla _{X}Y,Z)+g(Y,\nabla _{X}^{\ast }Z).
\]%
Thus $\nabla $ and $\nabla ^{\ast }$ are dual connections. Moreover, the
torsion tensor of $S(TM)$ with respect to $\nabla $ is given
\[
T^{\nabla }(X,Y)=\nabla _{X}Y-\nabla _{Y}X-[X,Y].
\]%
Using (\ref{d}) in the last equation we obtain $T^{\nabla }=0$. Similarly,
the torsion tensor of $S(TM)$ with respect to $\nabla ^{\ast }$ is equal to
zero. $\blacksquare $

\begin{prop}
\label{pr88} Let $(M,g)$ be a lightlike hypersurface of a statistical
manifold $(\widetilde{M},\widetilde{g},\widetilde{D},\widetilde{D}^{\ast })$%
. Then the following assertions are equivalent\/{\rm :}

\begin{enumerate}
\item[{\bf (i)}] The screen distribution $S(TM)$ is parallel.

\item[{\bf (ii)}] $C(X,Y)=0$ for all $X,Y\in\Gamma(S(TM))$.

\item[{\bf (iii)}] $C^{\ast }(X,Y)=0$ for all $X,Y\in\Gamma(S(TM))$.
\end{enumerate}
\end{prop}

\noindent {\bf Proof.} For any $X,Y\in\Gamma(S(TM))$, from Gauss-Weingarten
formulas and (\ref{18}), we obtain
\begin{equation}
g(D^{\ast }_{X}Y,N)=C^{\ast }(X,Y),
\end{equation}
\begin{equation}
g(D_{X}Y,N)=C(X,Y),
\end{equation}
Then, the proof is completed. $\blacksquare$

\begin{prop}
\label{pr8} Let $(M,g)$ be a lightlike hypersurface of a statistical
manifold $(\widetilde{M},\widetilde{g},\widetilde{D},\widetilde{D}^{\ast })$%
. Then the following assertions are equivalent\/{\rm :}

\begin{enumerate}
\item[{\bf (i)}] The screen distribution $S(TM)$ is integrable.

\item[{\bf (ii)}] $C(Y,X)=C(X,Y)$ for all $X,Y\in\Gamma(S(TM))$.

\item[{\bf (iii)}] $C^{\ast }(X,Y)=C^{\ast }(Y,X)$ for all $%
X,Y\in\Gamma(S(TM))$.
\end{enumerate}
\end{prop}

\noindent {\bf Proof.} For any $X,Y\in\Gamma(S(TM))$, from Gauss-Weingarten
formulas and (\ref{18}), we obtain
\begin{equation}
g([X,Y],N)=C(X,Y)-C(Y,X).
\end{equation}
\begin{equation}
g([X,Y],N)=C^{\ast }(X,Y)-C^{\ast }(Y,X).
\end{equation}
These equations prove our assertions. $\blacksquare$

\begin{defn-new}
{\rm (\cite{Furuhata-Has-2016}, \cite{Kurose-2002})} Let $(M,g)$ be a
hypersurface of a statistical manifold $(\widetilde{M},\widetilde{g},%
\widetilde{D},\widetilde{D}^{\ast })$.

\begin{enumerate}
\item[{\bf (i)}] $M$ is called totally geodesic with respect to $\widetilde{D%
}$ if $B=0$.

\item[{\bf (ii)}] $M$ is called totally geodesic with respect to $\widetilde{%
D}^{\ast }$ if $B^{\ast }=0$.

\item[{\bf (iii)}] $M$ is called totally tangentially umbilical with respect
to $\widetilde{D}$ if $B(X,Y)=kg(X,Y)$ for all $X,Y\in\Gamma(TM)$, where $k$
is smooth function.

\item[{\bf (iv)}] $M$ is called totally tangentially umbilical with respect
to $\widetilde{D}^{\ast }$ if $B^{\ast }(X,Y)=k^{\ast }g(X,Y)$, for any $%
X,Y\in\Gamma(TM)$, where $k^{\ast }$ is smooth function.

\item[{\bf (v)}] $M$ is called totally normally umbilical with respect to $%
\widetilde{D}$ if $A^{\ast }_{N}X=kX$ for any $X,Y\in\Gamma(TM)$, where $k$
is smooth function.

\item[{\bf (vi)}] $M$ is called totally normally umbilical with respect to $%
\widetilde{D}^{\ast }$ if $A_{N}X=k^{\ast }X$ for all $X,Y\in\Gamma(TM)$,
where $k^{\ast }$ is smooth function.
\end{enumerate}
\end{defn-new}

In view of (\ref{31}), (\ref{33}), (\ref{7}) and (\ref{8}), we have the
following proposition.

\begin{prop}
Let $(M,g)$ be a lightlike hypersurface of a statistical manifold $(%
\widetilde{M},\widetilde{g},\widetilde{D},\widetilde{D}^{\ast })$. Then the
following assertions are equivalent\/{\rm :}

\begin{enumerate}
\item[{\bf (i)}] $M$ is totally geodesic with respect to $\widetilde{D}$
(resp. $M$ is totally geodesic with respect to $\widetilde{D}^{\ast }$).

\item[{\bf (ii)}] $\overline{A}^{\ast }_{\xi}$ vanishes on $M$ (resp. $%
\overline{A}_{\xi}$ vanishes on $M$).

\item[{\bf (iii)}] $RadTM$ is a parallel distribution with respect to $%
\widetilde{D} $ (resp. $RadTM$ is a parallel distribution with respect to $%
\widetilde{D}^{\ast }$).

\item[{\bf (iv)}] $B^{\ast }(X,Y) = g(\overline{A}_{\xi}X,Y)$ (resp. $%
B(X,Y)=g(\overline{A}^{\ast }_{\xi}X,Y)$), for all $X,Y\in\Gamma(TM)$.
\end{enumerate}
\end{prop}

Next, we have the following

\begin{prop}
Let $(M,g)$ be a lightlike hypersurface of a statistical manifold $(%
\widetilde{M},\widetilde{g},\widetilde{D},\widetilde{D}^{\ast })$. Then the
following assertions are equivalent\/{\rm :}

\begin{enumerate}
\item[{\bf (i)}] $M$ is totally geodesic with respect to $\widetilde{D}$ and
$\widetilde{D}^{\ast }$.

\item[{\bf (ii)}] $\overline{A}_{\xi}X=\overline{A}_{\xi}^{\ast }X=0$ for
all $X\in\Gamma(TM)$.

\item[{\bf (iii)}] $D_{X}g + D^{\ast }_{X}g=0$ for all $X\in\Gamma(TM)$.

\item[{\bf (iv)}] $D_{X}\xi+D^{\ast }_{X}\xi\in\Gamma(RadTM)$ for all $%
X\in\Gamma(TM)$.
\end{enumerate}
\end{prop}

\noindent {\bf Proof.} From (\ref{34}), (\ref{7}) and (\ref{8}) we get the
equivalence of {\bf (i)} and {\bf (ii)}. The equation (\ref{th-LH-SM-2})
implies the equivalence of {\bf (i)} and {\bf (iii)}. Next, by using (\ref%
{31}) and (\ref{33}) we have the equivalence of {\bf (ii)} and {\bf (iv)}. $%
\blacksquare$

\begin{th}
Let $(M,g)$ be a lightlike hypersurface of a statistical manifold $(%
\widetilde{M},\widetilde{g},\widetilde{D},\widetilde{D}^{\ast })$. Then, $M$
is totally tangentially umbilical with respect to $\widetilde{D}$ and $%
\widetilde{D}^{\ast }$ if and only if
\[
\overline{A}^{\ast }_{\xi}X+\overline{A}_{\xi}X=\rho X,\;\forall
X\in\Gamma(TM),
\]
where $\rho$ is smooth function.
\end{th}

\noindent {\bf Proof.} Using (\ref{7}) and (\ref{8}) we obtain
\begin{eqnarray}
kg(X,Y)=g(\overline{A}^{\ast }_{\xi}X,Y)+B^{\ast }(X,\xi),  \label{t11}
\end{eqnarray}
and
\begin{eqnarray}
k^{\ast }g(X,Y)=g(\overline{A}_{\xi}X,Y)+B(X,\xi).  \label{t222}
\end{eqnarray}
If we add the equations (\ref{t11}) and (\ref{t222}) side by side and using (%
\ref{34}) we complete the proof. $\blacksquare$

\begin{prop}
Let $(M,g)$ be a lightlike hypersurface of a statistical manifold $(%
\widetilde{M},\widetilde{g},\widetilde{D},\widetilde{D}^{\ast })$. If $M$ is
totally normally umbilical with respect to $\widetilde{D}$ and $\widetilde{D}%
^{\ast }$. Then
\[
C(X,PY)+C^{\ast }(X,PY)=0,\;\forall X\in\Gamma(TM).
\]
\end{prop}

\noindent {\bf Proof.} Let $k$ and $k^{\ast }$ be smooth functions and let $%
A_{N}^{\ast }X=kX$ and $A_{N}X=k^{\ast }X$, then using (\ref{34}) we get $%
k+k^{\ast }=0$. Thus, from (\ref{18}) proof is completed. $\blacksquare$

\medskip

It is known that $M$ is screen locally conformal lightlike hypersurface of a
statistical manifold $\widetilde{M}$ if
\begin{eqnarray}
A_{N} = \varphi \overline{A}_{\xi}^{\ast },\;A^{\ast }_{N} = \varphi^{\ast }
\overline{A}_{\xi},  \label{sc}
\end{eqnarray}
where $\varphi$ and $\varphi^{\ast }$ are non-vanishing smooth functions on $%
M$. Using (\ref{18}) and (\ref{sc}) we get the following proposition.

\begin{prop}
Let $(M,g)$ be a lightlike hypersurface of a statistical manifold $(%
\widetilde{M},\widetilde{g},\widetilde{D},\widetilde{D}^{\ast })$. Then, $M$
is screen locally conformal if and only if
\[
C(X,Y)+C^{\ast }(X,Y)=\sigma(B(X,Y)+B^{\ast }(X,Y)),\;\forall X,Y\in\Gamma
(S(TM)),
\]
where $\sigma$ is non-vanishing smooth functions on $M$.
\end{prop}

Now, we give an example.

\begin{ex-new}
Let $(R_{2}^{4},\widetilde{g})$ be a $4$-dimensional semi-Euclidean space
with signature \allowbreak $(-,-,+,+)$ of the canonical basis $%
(\partial_{0},\ldots ,\partial_{3})$. Consider a hypersurface $M$ of $%
R_{2}^{4}$ given by
\[
x_{0}=x_{1}+\sqrt{2}\sqrt{x_{2}^{2}+x_{3}^{2}}.
\]
For simplicity, we set $f=\sqrt{x_{2}^{2}+x_{3}^{2}}$. It is easy to check
that $M$ is a lightlike hypersurface whose radical distribution $RadTM$ is
spanned by
\[
\xi=f(\partial_{0}-\partial_{1})+\sqrt{2}(x_{2}\partial_{2}+x_{3}%
\partial_{3}).
\]
Then the lightlike transversal vector bundle is given by
\[
ltr(TM)=Span\{{N=\frac{1}{4f^{2}}\{f(-\partial_{0}+\partial_{1})+\sqrt{2}%
(x_{2}\partial_{2}+x_{3}\partial_{3})\}}\}.
\]
It follows that the corresponding screen distribution $S(TM)$ is spanned by
\[
\{W_{1}=\partial_{0}+\partial_{1},\;W_{2}=-x_{3}\partial_{2}+x_{2}%
\partial_{3}\}.
\]
Then, by direct calculations we obtain
\[
\widetilde{\nabla}_{X}W_{1}=\widetilde{\nabla}_{W_{1}}X=0,
\]
\[
\widetilde{\nabla}_{W_{2}}W_{2}=-x_{2}\partial_{2}-x_{3}\partial_{3},
\]
\[
\widetilde{\nabla}_{\xi}\xi=\sqrt{2}\xi,\;\widetilde{\nabla}_{W_{2}}\xi=%
\widetilde{\nabla}_{\xi}W_{2}=\sqrt{2}W_{2},
\]
for any $X\in\Gamma(TM)$ \cite{Duggal-Sahin-2010}.

We define an affine connection $\widetilde{D}$ as follows
\begin{eqnarray}
&&\widetilde{D}_{X}W_{1}=\widetilde{D}_{W_{1}}X=0,\;\widetilde{D}%
_{W_{2}}W_{2}=-2x_{2}\partial_{2}  \nonumber \\
&&\widetilde{D}_{\xi}\xi=\sqrt{2}\xi-\sqrt{2}N,  \label{denk1} \\
&& \widetilde{D}_{W_{2}}\xi=\widetilde{D}_{\xi}W_{2}=\sqrt{2}W_{2}-\sqrt{2}%
W_{1}.  \nonumber
\end{eqnarray}
Then using (\ref{le}) we obtain
\begin{eqnarray}
&&\widetilde{D}^{\ast }_{X}W_{1}=\widetilde{D}^{\ast }_{W_{1}}X=0,\;%
\widetilde{D}^{\ast }_{W_{2}}W_{2}=-2x_{3}\partial_{3}  \nonumber \\
&&\widetilde{D}^{\ast }_{\xi}\xi=\sqrt{2}\xi+\sqrt{2}N,  \label{denk2} \\
&&\widetilde{D}^{\ast }_{W_{2}}\xi=\widetilde{D}^{\ast }_{\xi}W_{2}=\sqrt{2}%
W_{2}+\sqrt{2}W_{1}.  \nonumber
\end{eqnarray}
Then $(R_{2}^{4},\widetilde{g},\widetilde{D},\widetilde{D}^{\ast })$ is a
statistical manifold. Thus, by using Gauss formulas (\ref{eq-D-tilde-X-Y})
and (\ref{eq-D-tilde-star-X-Y}) we obtain
\begin{eqnarray}
&&B(X,W_{1})=B(W_{1},X)=0,  \nonumber \\
&&B(W_{2},W_{2})=-2\sqrt{2}x_{2}^{2},\; B(\xi,\xi)=-\sqrt{2}  \label{denk3}
\\
&&B(X,W_{2})=B(W_{2},X)=0,  \nonumber
\end{eqnarray}
and
\begin{eqnarray}
&&B^{\ast }(X,W_{1})=B^{\ast }(W_{1},X)=0,  \nonumber \\
&&B^{\ast }(W_{2},W_{2})=-2\sqrt{2}x_{3}^{2},\; B^{\ast }(\xi,\xi)=\sqrt{2}
\label{denk4} \\
&&B^{\ast }(X,W_{2})=B^{\ast }(W_{2},X)=0.  \nonumber
\end{eqnarray}
The equations (\ref{denk1}), (\ref{denk2}), (\ref{denk3}) and (\ref{denk4})
imply that induced connections $D$ and $D^{\ast }$ are symmetric connections
and the second fundamental forms $B$ and $B^{\ast }$ are symmetric. This
proves Proposition~\ref{pro5}. Moreover, the equations $B(\xi,\xi)=-\sqrt{2}$
and $B^{\ast }(\xi,\xi)=\sqrt{2}$ show the accuracy of the Proposition~\ref%
{pr6}.

Using (\ref{denk1}), (\ref{denk2}), (\ref{denk3}) and (\ref{denk4}) we get
\begin{eqnarray}
&&D_{X}W_{1}=D_{W_{1}}X=0,\;D_{\xi}\xi=\sqrt{2}\xi,  \nonumber \\
&&D_{W_{2}}W_{2}=\frac{\sqrt{2}x_{2}^{2}}{2f}(-\partial_{0}+\partial_{1})+%
\frac{1}{4f^{2}}\{(4x_{2}^{3}-2x_{2})\partial_{2}+4x_{3}x_{2}^{2}%
\partial_{3})\},  \label{denk5} \\
&&D_{W_{2}}\xi=D_{\xi}W_{2}=\sqrt{2}W_{2}-\sqrt{2}W_{1},  \nonumber
\end{eqnarray}
and
\begin{eqnarray}
&&D^{\ast }_{X}W_{1}=D^{\ast }_{W_{1}}X=0,\;D^{\ast }_{\xi}\xi=\sqrt{2}\xi,
\nonumber \\
&&D^{\ast }_{W_{2}}W_{2}=\frac{\sqrt{2}x_{3}^{2}}{2f}(-\partial_{0}+%
\partial_{1})+\frac{1}{4f^{2}}\{4x_{3}^{2}x_{2}%
\partial_{2}+(4x_{3}^{3}-2x_{3})\partial_{3})\},  \label{denk6} \\
&&D^{\ast }_{W_{2}}\xi=D^{\ast }_{\xi}W_{2}=\sqrt{2}W_{2}+\sqrt{2}W_{1}.
\nonumber
\end{eqnarray}
In the equation (\ref{eq-dual-con}), if we choose $X=W_{2}$, $Y=W_{2}$ and $%
Z=\xi$, (\ref{denk5}) and (\ref{denk6}) indicate that induced connections $%
D^{\ast }$ and $D$ are not dual connections. This verifies Theorem~\ref%
{th-LH-SM-1}.

From (\ref{c}) and (\ref{d}), we have
\begin{equation}
C(X,W_{1})=C(W_{1},X)=0,\;C(W_{2},W_{2})=-\frac{\sqrt{2}}{2}(\frac{x_{2}}{f}%
)^{2},\;C(\xi,W_{2})=0  \label{c11}
\end{equation}
and
\begin{eqnarray}
C^{\ast }(X,W_{1})=C^{\ast }(W_{1},X)=0,\;C^{\ast }(W_{2},W_{2}) = -\frac{%
\sqrt{2}}{2}(\frac{x_{3}}{f})^{2},\;C^{\ast }(\xi,W_{2})=0.  \label{c22}
\end{eqnarray}
From (\ref{c11}) and (\ref{c22}), we say that $C$ and $C^{\ast }$ are
symmetric. Thus we have Proposition~\ref{pr8}.

Using (\ref{denk5}) and (\ref{denk6}) in (\ref{c}) and (\ref{d}) we obtain
\begin{eqnarray}
&&\nabla_{X}W_{1}=\nabla_{W_{1}}X=0,  \nonumber \\
&&\nabla_{W_{2}}W_{2}=\frac{1}{f^{2}}\{(2x_{2}^{3}-\frac{x_{2}}{2}%
)\partial_{2}+2x_{3}x_{2}^{2}\partial_{3}\},  \label{denk56} \\
&&\nabla_{\xi}W_{2}=\sqrt{2}W_{2}-\sqrt{2}W_{1},  \nonumber
\end{eqnarray}
and
\begin{eqnarray}
&&\nabla^{\ast }_{X}W_{1}=\nabla^{\ast }_{W_{1}}X=0,  \nonumber \\
&&\nabla^{\ast }_{W_{2}}W_{2}=\frac{1}{f^{2}}\{2x_{3}^{2}x_{2}%
\partial_{2}+(2x_{3}^{3}-\frac{x_{3}}{2})\partial_{3}\},  \label{denk66} \\
&&\nabla^{\ast }_{\xi}W_{2}=\sqrt{2}W_{2}+\sqrt{2}W_{1}.  \nonumber
\end{eqnarray}
From (\ref{denk56}) and (\ref{denk66}), the torsion tensors vanish with
respect to $\nabla$ and $\nabla^{\ast }$. Furthermore, this equations
provides (\ref{eq-dual-con}). Thus, $\nabla$ and $\nabla^{\ast }$ are dual
connections. This situation verifies Proposition~\ref{pr7}.
\end{ex-new}

\section{Curvature tensors of a lightlike hypersurface of a statistical
manifold\label{sect-curvature}}

We denote by $\widetilde{R}$ and $\widetilde{R}^{\ast }$ the curvature
tensor of $\widetilde{D}$ and $\widetilde{D}^{\ast }$, respectively. The
curvature tensors satisfy
\[
\widetilde{g}(\widetilde{R}^{\ast }(X,Y)Z,W)=-\widetilde{g}(\widetilde{R}%
(X,Y)W,Z).
\]%
Using Gauss-Weingarten formulas, the curvature tensors $\widetilde{R}$ and $%
\widetilde{R}^{\ast }$ of the connection $\widetilde{D}$ and $\widetilde{D}%
^{\ast }$ are given by
\begin{eqnarray}
\widetilde{R}(X,Y)Z &=&R(X,Y)Z-B(Y,Z)A_{N}^{\ast }X+B(X,Z)A_{N}^{\ast }Y
\nonumber \\
&&+(B(Y,Z)\tau ^{\ast }(X)-B(X,Z)\tau ^{\ast }(Y))N  \nonumber \\
&&+((D_{X}B)(Y,Z)-(D_{Y}B)(X,Z))N,  \label{R}
\end{eqnarray}%
and
\begin{eqnarray}
\widetilde{R}^{\ast }(X,Y)Z &=&R^{\ast }(X,Y)Z-B^{\ast }(Y,Z)A_{N}X+B^{\ast
}(X,Z)A_{N}Y  \nonumber \\
&+&(B^{\ast }(Y,Z)\tau (X)-B^{\ast }(X,Z)\tau (Y))N  \nonumber \\
&+&((D_{X}^{\ast }B^{\ast })(Y,Z)-(D_{Y}^{\ast }B^{\ast })(X,Z))N,
\label{R*}
\end{eqnarray}%
where $R$ and $R^{\ast }$ are the curvature tensor with respect to $D$ and $%
D^{\ast }$, respectively. Consider curvature tensors $\widetilde{R}$ and $%
\widetilde{R}^{\ast }$ of type $(0,4)$. From the above equation and the
Gauss-Weingarten equations for $M$ and $S(TM)$ we obtain
\begin{eqnarray}
g(\widetilde{R}(X,Y)Z,PW) &=&g(R(X,Y)Z,PW)-B(Y,Z)C^{\ast }(X,PW)  \nonumber
\\
&&+B(X,Z)C^{\ast }(Y,PW),  \label{g1}
\end{eqnarray}%
\begin{eqnarray}
g(\widetilde{R}^{\ast }(X,Y)Z,PW) &=&g(R^{\ast }(X,Y)Z,PW)-B^{\ast
}(Y,Z)C(X,PW)  \nonumber \\
&&+B^{\ast }(X,Z)C(Y,PW),  \label{g2}
\end{eqnarray}%
\begin{eqnarray}
g(\widetilde{R}(X,Y)Z,\xi ) &=&B(Y,Z)\tau ^{\ast }(X)-B(X,Z)\tau ^{\ast }(Y)
\nonumber \\
&&+(D_{X}B)(Y,Z)-(D_{Y}B)(X,Z),  \label{g3}
\end{eqnarray}%
\begin{eqnarray}
g(\widetilde{R}^{\ast }(X,Y)Z,\xi ) &=&B^{\ast }(Y,Z)\tau (X)-B^{\ast
}(X,Z)\tau (Y)  \nonumber \\
&&+(D_{X}^{\ast }B^{\ast })(Y,Z)-(D_{Y}^{\ast }B^{\ast })(X,Z),  \label{g4}
\end{eqnarray}%
\begin{eqnarray}
g(\widetilde{R}(X,Y)Z,N) &=&g(R(X,Y)Z,N)-B(Y,Z)g(A_{N}^{\ast }X,N)  \nonumber
\\
&&+B(X,Z)g(A_{N}^{\ast }Y,N),  \label{g5}
\end{eqnarray}%
\begin{eqnarray}
g(\widetilde{R}^{\ast }(X,Y)Z,N) &=&g(R^{\ast }(X,Y)Z,N)-B^{\ast
}(Y,Z)g(A_{N}X,N)  \nonumber \\
&&+B^{\ast }(X,Z)g(A_{N}Y,N),  \label{g6}
\end{eqnarray}%
\begin{eqnarray}
g(\widetilde{R}(X,Y)\xi ,N) &=&g(R(X,Y)\xi ,N)-B(Y,\xi )g(A_{N}^{\ast }X,N)
\nonumber \\
&&+B(X,\xi )g(A_{N}^{\ast }Y,N),  \label{g7}
\end{eqnarray}%
\begin{eqnarray}
g(\widetilde{R}^{\ast }(X,Y)\xi ,N) &=&g(R^{\ast }(X,Y)\xi ,N)-B^{\ast
}(Y,\xi )g(A_{N}X,N)  \nonumber \\
&&+B^{\ast }(X,\xi )g(A_{N}Y,N),  \label{g8}
\end{eqnarray}%
where
\[
g(R(X,Y)\xi ,N)=C(Y,\overline{A}_{\xi }X)-C(X,\overline{A}_{\xi }Y)-2d\tau
(X,Y),
\]%
\[
g(R^{\ast }(X,Y)\xi ,N)=C^{\ast }(Y,\overline{A}_{\xi }^{\ast }X)-C^{\ast
}(X,\overline{A}_{\xi }^{\ast }Y)-2d\tau (X,Y).
\]

Now, let $M$ be a lightlike hypersurface of a $(m+2)$-dimensional
statistical manifold $\widetilde{M}$. We consider the local
quasi-orthonormal basis $\{E_{i},\xi,N\},\ i=1,\ldots m,$ of
$\widetilde{M}$ along $M$, where $\{E_{1},\ldots ,E_{m}\}$ is an
orthonormal basis of $\Gamma(S(TM))$. Then, we obtain
\begin{equation}
R^{D(0,2)}(X,Y) =\sum_{i=1}^{m}\varepsilon_{i}g(R(X,E_{i})Y,E_{i})+
\widetilde{g}(R(X,\xi)Y,N),  \label{ri}
\end{equation}
where $\varepsilon_{i}$ denotes the causal character $(\mp1)$ of
respective vector field $E_{i}$. Using Gauss-Weingarten equations we
have
\begin{equation}
g(R(X,E_{i})Y,E_{i})=g(\widetilde{R}(X,E_{i})Y,E_{i})+B(E_{i},Y)C^{\ast
}(X,E_{i})-B(X,Y)C^{\ast }(E_{i},E_{i})
\end{equation}
Substituting this in (\ref{ri}), using (\ref{18}) and (\ref{7}) we
obtain
\begin{equation}
R^{D(0,2)}(X,Y)=\widetilde{Ric}(X,Y)-B(X,Y)trA_{N}^{\ast }+g(A_{N}^{\ast }X,%
\overline{A}_{\xi}^{\ast }Y)+g(R(X,\xi)Y,N)  \label{ric}
\end{equation}
where $\widetilde{Ric}(X,Y)$ is the Ricci tensor of $\widetilde{M}$
with respect to $\widetilde{D}$. Similarly, dual tensor of $M$ with
respect to $D^{\ast }$ as follows:
\begin{equation}
R^{D^{\ast }(0,2)}(X,Y)=\widetilde{Ric}^{\ast }(X,Y)-B^{\ast
}(X,Y)trA_{N}+g(A_{N}X,\overline{A}_{\xi}Y)+g(R^{\ast }(X,\xi)Y,N)
\label{ric*}
\end{equation}
From First Bianchi identities and (\ref{ric}) we get
\begin{eqnarray}
R^{D(0,2)}(X,Y)-R^{D(0,2)}(Y,X)&=&\sum_{i=1}^{m}%
\varepsilon_{i}((B(E_{i},Y)C^{\ast }(X,E_{i})-B(E_{i},X)C^{\ast
}(Y,E_{i}) \nonumber \\
&&+g(\widetilde{R}(X,Y)E_{i},E_{i})) +g(\widetilde{R}(X,Y)\xi,N).
\label{ricf}
\end{eqnarray}
Therefore, $R^{D(0,2)}$ is not symmetric.

The statistical manifold $(\widetilde{M},\widetilde{g})$ is called
of constant curvature $c$ if
\begin{equation}
\widetilde{R}(X,Y)Z=c(Y,Z)X-g(X,Z)Y .  \label{s}
\end{equation}
Moreover, if $(\widetilde{D}, \widetilde{g})$ is a statistical
structure of constant c, then $(\widetilde{D}^{\ast },
\widetilde{g})$ is also a statistical structure of constant c
\cite{Furuhata-2011}. Then, using (\ref{18}), (\ref{7}),
({\ref{g7}}) and (\ref{s}) in (\ref{ricf}) we have
\begin{equation}
R^{D(0,2)}(X,Y)-R^{D(0,2)}(Y,X)=C^{\ast }(X,\overline{A}_{\xi}^{\ast
}Y)-C^{\ast }(Y,\overline{A}_{\xi}^{\ast }X),  \label{rics}
\end{equation}
and similarly
\begin{eqnarray}
R^{D^{\ast }(0,2)}(X,Y)-R^{D^{\ast }(0,2)}(Y,X)=C(X,\overline{A}_{\xi}Y)-C(Y,%
\overline{A}_{\xi}X).  \label{ricss}
\end{eqnarray}
Then we have the following theorem

\begin{th}
Let $(M,g)$ be a lightlike hypersurface of a statistical manifold
$(\widetilde{M}^{n+2}(c),\widetilde{g})$ of constant sectional
curvature c. Then the following assertions are true\/{\rm :}
\begin{enumerate}
\item[{\bf (i)}] The tensor $R^{D(0,2)}(X,Y)$ is symmetric if and
only if
\[
C^{\ast }(X,\overline{A}_{\xi}^{\ast }Y)=C^{\ast
}(Y,\overline{A}_{\xi}^{\ast }X).
\]

\item[{\bf (ii)}]The tensor $R^{D^{\ast }(0,2)}(X,Y)$ is symmetric if and
only if
\[
C(X,\overline{A}_{\xi}Y)=C(Y,\overline{A}_{\xi}X).
\]
\end{enumerate}
\end{th}

Thus, in view of Propositon~\ref{pr88}, we have the following:

\begin{cor}
Let $(M,g)$ be a lightlike hypersurface of a statistical manifold
$(\widetilde{M}^{n+2}(c),\widetilde{g})$ of constant sectional
curvature c. If $S(TM)$ is parallel then the tensor $R^{D(0,2)}$ and
$R^{D^{\ast }(0,2)}$ are symmetric with respect to connections $D$
and $D^{\ast }$, respectively.
\end{cor}

Now, let $M$ be a $(m + 1)$-dimensional lightlike hypersurface of a
Lorentzian space form $\widetilde{M}(c)$. We know that degenerate
scalar curvature $\sigma$ is given by
\begin{equation}
\sigma(u)=r_{S(TM)}+\sum_{i=1}^{m}\{\kappa_{i}^{null}+\kappa_{iN}\},
\label{scal}
\end{equation}
where $\kappa_{iN}=\widetilde{g}(R(E_{i},\xi)E_{i},N)$ and $
r_{S(TM)}=\sum_{i,j=1}^{m}\kappa _{ij}$ is called screen scalar
curvature \cite{Gulbahar-KK-2013}.

Using (\ref{g1}) and (\ref{ri}) we have
\begin{eqnarray}
g(R(X,Y)X,PW)&=&c(g(Y,Z)g(X,PW)-g(X,Z)g(Y,PW))  \nonumber \\
&+&B(Y,Z)C^{\ast }(X,PW)-B(X,Z)C^{\ast }(Y,PW).  \label{osu}
\end{eqnarray}
Then (\ref{scal}) and $(\ref{osu})$ gives the following equation
\begin{eqnarray}
r_{S(TM)}&=&\sum_{i,j=1}^{m}c\{{%
g(E_{j},E_{i})g(E_{i},E_{j})-g(E_{i},E_{i})g(E_{j},E_{j})}\}  \nonumber \\
&+&B(E_{j},E_{i})C^{\ast }(E_{i},E_{j})-B(E_{i},E_{i})C^{\ast }(E_{j},E_{j})
\nonumber \\
&=&cm(1-m)+\sum_{i,j}^{m}(B_{ji}C^{\ast }_{ij}-B_{ii}C^{\ast }_{jj}),
\label{47*}
\end{eqnarray}
in where $B_{ji}=B(E_{j},E_{i})$ and $C^{\ast }_{ij}=C^{\ast
}(E_{i},E_{j})$. If we replace $X=\xi=Y$ in (\ref{ri}) and using
(\ref{g5}), we obtain
\begin{equation}
R^{D(0,2)}(\xi,\xi) =\kappa_{i}^{null}.  \label{12*}
\end{equation}
Thus, (\ref{osu}) is gives us
\begin{equation}
\sum_{i=1}^{m}\kappa_{i}^{null}=\sum_{i=1}^{m}\{B(E_{i},\xi)C^{\ast
}(\xi,E_{i})-B(\xi,\xi)C^{\ast }(E_{i},E_{i})\}.  \label{nul}
\end{equation}
Moreover, using (\ref{g5}), we get
\begin{equation}
\sum_{i=1}^{m}\kappa_{iN}=-cm-\sum_{i=1}^{m}B(\xi,E_{i})g(A^{\ast
}_{N}E_{i},N)-B(E_{i},E_{i})g(A^{\ast }_{N}\xi,N).  \label{k1n}
\end{equation}


\noindent O\u{g}uzhan Bahad\i r

\noindent Department of Mathematics,

\noindent Faculty of Arts and Sciences,

\noindent Kahramanmaras Sut\c{c}u Imam University

\noindent Kahramanmaras, T\"{u}rkiye

\noindent Email: oguzbaha@gmail.com

\medskip

\noindent Mukut Mani Tripathi

\noindent Department of Mathematics,

\noindent Institute of Sciences,

\noindent Banaras Hindu University

\noindent Varanasi 221005, India

\noindent Email: mmtripathi66@yahoo.com

\end{document}